\documentclass[leqno,12pt]{article}%
\usepackage{amssymb}
\usepackage{amsmath}
\usepackage{amsfonts}
\usepackage{graphicx}%
\providecommand{\U}[1]{\protect\rule{.1in}{.1in}}
\begin{document}

\begin{center}
{\large MOD-2 EQUIVALENCE OF THE \textit{K}-THEORETIC EULER AND SIGNATURE
CLASSES}{\Large \textbf{\footnotetext{2000 Mathematics Subject
Classification:\ 19L99, 58J05, 19K56, 57R20.}}}

\bigskip

\textsc{James F. Davis and Pisheng Ding}

\bigskip
\end{center}

\begin{quote}
{\small This note proves that, as }$K${\small -theory elements, the symbol
classes of the de Rham operator and the signature operator on a closed
manifold of even dimension\ are congruent mod 2. An equivariant generalization
is given pertaining to the equivariant Euler characteristic and the
multi-signature.}
\end{quote}

\bigskip

\begin{center}
\textsc{1. Introduction}

\medskip
\end{center}

It is well-known that the Euler characteristic $\chi(M)$ and the signature
$\operatorname*{Sign}(M)$ of a closed oriented manifold $M$ of dimension $4n$
are two integers of the same parity. This fact is an easy consequence of
Poincar\'{e} duality and we briefly recall its proof. Let $\beta_{k}=\dim
H_{k}(M;{\mathbb{R}})$. By Poincar\'{e} duality, $\beta_{k}=\beta_{4n-k}$ and
there is a non-degenerate bilinear form (the \textquotedblleft intersection
form\textquotedblright) on $H_{2n}(M;{\mathbb{R}})$. Let $\beta_{2n}^{+}$
(respectively $\beta_{2n}^{-}$) be the dimension of a maximal subspace of
$H_{2n}(M;{\mathbb{R}})$ on which the form is positive definite (respectively
negative definite). By definition, $\operatorname*{Sign}(M)=\beta_{2n}%
^{+}-\beta_{2n}^{-}$. The following mod-2 congruence relation follows:%
\[
\chi(M)=\sum_{k=0}^{4n}(-1)^{k}\beta_{k}\equiv\beta_{2n}^{+}+\beta_{2n}%
^{-}+\sum_{k=0}^{2n-1}2\beta_{k}\equiv\beta_{2n}^{+}-\beta_{2n}^{-}%
=\operatorname*{Sign}(M)\text{.}%
\]

We prove an analogous result on the level of the symbols (as $K$-theory
classes) of the de Rham and signature operators on even-dimensional manifolds;
we also give an equivariant generalization.

We first recall the construction of the two operators (cf. \cite{AS3}).

Let $M$ be a closed oriented smooth manifold of dimension $2n$. Equip $M$ with
a Riemannian metric. (The symbols, as $K$-theory classes, of the two operators
will be independent of the choice of Riemannian metric and hence are smooth
invariants. This is due to the fact that, for any metrics $g_{0}$ and $g_{1}$
and for $t\in\lbrack0,1]$, $tg_{0}+(1-t)g_{1}$ is a Riemannian metric.) Let%
\[
\Omega^{\ast}=\Gamma(\Lambda^{\ast}(T^{\ast}M\otimes\mathbb{C}))\text{,}%
\]
the space of smooth sections of the exterior algebra bundle $\Lambda^{\ast
}=\Lambda^{\ast}(T_{\mathbb{C}}^{\ast}M)$ associated with the complexified
cotangent bundle $T_{\mathbb{C}}^{\ast}M=T^{\ast}M\otimes\mathbb{C}$; i.e.,
$\Omega^{\ast}$ is the space of \textit{complex} differential forms. Let%
\[
D=d+\delta:\Omega^{\ast}\rightarrow\Omega^{\ast}%
\]
where $d$ is the exterior derivative and $\delta$ its adjoint with respect to
the Hermitian product on $\Omega^{\ast}$ induced by the Riemannian metric on
$M$. The de Rham operator $D^{0}$ is defined to be the restriction of $D$ to
the subspace $\Omega^{\text{even}}$ of even-degree forms and thus takes value
in the subspace $\Omega^{\text{odd}}$ of odd-degree forms:%
\[
D^{0}=D|_{\Omega^{\text{even}}}:\Omega^{\text{even}}=\Gamma(\Lambda
^{\text{even}})\rightarrow\Omega^{\text{odd}}=\Gamma(\Lambda^{\text{odd}%
})\text{.}%
\]

The bundle map Hodge star $\ast:\Lambda^{k}\rightarrow\Lambda^{2n-k}$ is
defined on each fiber $(\Lambda^{k})_{x}$ by%
\[
\alpha\wedge\ast\beta=<\alpha,\beta>\cdot\operatorname*{vol}(M)_{x}\text{\quad
for }\alpha,\beta\in(\Lambda^{k})_{x}%
\]
where $\operatorname*{vol}(M)\in\Omega^{2n}$ is the volume form. The Hodge
star has the property that%
\[
\ast(\ast\alpha)=(-1)^{k}\alpha\text{\quad for }\alpha\in\Lambda^{k}\text{.}%
\]
Define $\tau:\Lambda^{k}\rightarrow\Lambda^{2n-k}$ by letting%
\[
\tau=i^{n+k\left(  k-1\right)  }\ast\text{.}%
\]
It is easy to verify that $\tau$ is an involution on $\Lambda^{\ast}$. Then
$\tau$ decomposes $\Lambda^{\ast}$ into $\Lambda^{+}\oplus\Lambda^{-}$, the
$+1$ and $-1$ eigenbundles. The map $\tau$ induces an involution on
$\Omega^{\ast}$ (which we also call $\tau$) and decomposes $\Omega^{\ast}$
into $\Omega^{+}\oplus\Omega^{-}$, where $\Omega^{\pm}$ are the ($\pm
1$)-eigenspaces. Note that $\Omega^{\pm}=\Gamma(\Lambda^{\pm})$ and that $D$
interchanges $\Omega^{+}$ and $\Omega^{-}$ (since $D\tau=-\tau D$). The
signature operator $D^{+}$ is defined to be the restriction of $D$ to
$\Omega^{+}$:%
\[
D^{+}=D|_{\Omega^{+}}:\Omega^{+}=\Gamma(\Lambda^{+})\rightarrow\Omega
^{-}=\Gamma(\Lambda^{-})\text{.}%
\]

We now recall the $K$-theoretic symbol class associated with an elliptic
differential operator. Suppose $E_{1}$ and $E_{2}$ are two complex vector
bundles over $M$. Let $\pi:T_{\mathbb{C}}^{\ast}M\rightarrow M$ be the bundle
projection. Associated with a differential operator $P:\Gamma(E_{1}%
)\rightarrow\Gamma(E_{2})$, there is the (leading) symbol of $P$,%
\[
\sigma(P):\pi^{\ast}E_{1}\rightarrow\pi^{\ast}E_{2}\text{.}%
\]
$\sigma(P)$ is a bundle homomorphism. If, over the complement of the zero
section of $T_{\mathbb{C}}^{\ast}M$, $\sigma(P)$ is a bundle isomorphism, $P$
is then said to be \textit{elliptic}. The symbol $\sigma(P)$ of an elliptic
operator $P$ determines a class $\left[  \sigma(P)\right]  $, the
$K$\textit{-theoretic symbol class of }$P$, in%
\[
K(T_{\mathbb{C}}^{\ast}M)=\widetilde{K}(D(T_{\mathbb{C}}^{\ast}%
M)/S(T_{\mathbb{C}}^{\ast}M))\text{,}%
\]
where $D(T_{\mathbb{C}}^{\ast}M)$ and $S(T_{\mathbb{C}}^{\ast}M)$ denote the
closed unit disc bundle and the unit sphere bundle associated with
$T_{\mathbb{C}}^{\ast}M$. We will review in more detail the transition from
$\sigma(P)$ to $\left[  \sigma(P)\right]  $ in \S 2.

It is a standard fact that $D^{0}$ and $D^{+}$ are elliptic. We call
$[\sigma(D^{0})]$ the $K$\textit{-theoretic Euler class of }$M$ and
$[\sigma(D^{+})]$ the $K$\textit{-theoretic signature class of }$M$. It is
their relationship in the abelian group $K(T_{\mathbb{C}}^{\ast}M)$ that our
Main Theorem pertains to.

\medskip

\textsc{Main Theorem. }\textit{If }$\dim M$\textit{ is even, then }$\left[
\sigma\left(  D^{0}\right)  \right]  \equiv\left[  \sigma\left(  D^{+}\right)
\right]  \operatorname{mod}2K(T_{\mathbb{C}}^{\ast}M)$\textit{.}

\medskip

(Here, $2K(T_{\mathbb{C}}^{\ast}M)$ means $\left\{  \alpha+\alpha:\alpha\in
K(T_{\mathbb{C}}^{\ast}M)\right\}  $.)

\medskip

The Main Theorem implies, as a corollary, the mod-2 congruence between
$\chi(M)$ and $\operatorname*{Sign}(M)$. To see this, consider the index
homomorphism $\operatorname*{Index}:K(T_{\mathbb{C}}^{\ast}M)\rightarrow
\mathbb{Z}$ and recall that $\operatorname*{Index}(\left[  \sigma
(D^{0})\right]  )=\chi(M)$ and $\operatorname*{Index}(\left[  \sigma
(D^{+})\right]  )=\operatorname*{Sign}(M)$.

\bigskip

\begin{center}
\textsc{2. Preliminaries on \textit{K}-theory}

\medskip
\end{center}

We review and establish some $K$-theoretic results.

Let $\left(  X,A\right)  $ be a pair of connected compact Hausdorff spaces.
Let $E_{1}$ and $E_{2}$ be two complex vector bundles over $X$ with a bundle
isomorphism $\sigma:E_{1}|_{A}\rightarrow E_{2}|_{A}$. The triple $\left(
E_{1},E_{2};\sigma\right)  $ determines a class $\left[  E_{1},E_{2}%
;\sigma\right]  $ in $\widetilde{K}(X/A)$. We first recall its construction.

Let $X_{i}=X\times\left\{  i\right\}  $ and $A_{i}=A\times\left\{  i\right\}
$, $i=1,2$; let $Y=X_{1}\cup_{g}X_{2}$, with $g:\left(  a,1\right)
\mapsto\left(  a,2\right)  $ for $a\in A$. (For notational convenience, we
regard $E_{i}$ as a bundle not only over $X$ but also over $X_{1}$ and $X_{2}%
$.) We first construct bundles $E_{i,j}$'s over $Y$. To produce $E_{i,j}$, we
glue $E_{i}|_{X_{1}}$ and $E_{j}|_{X_{2}}$ via $\varepsilon_{i,j}%
:E_{i}|_{A_{1}}\rightarrow E_{j}|_{A_{2}}$ where $\varepsilon_{1,2}=\sigma$,
$\varepsilon_{2,1}=\sigma^{-1}$, and $\varepsilon_{i,i}=\operatorname{Id}%
_{E_{i}|_{A}}$.

Consider $E_{1,2}-E_{2,2}\in\widetilde{K}(Y)$. Evidently, its restriction to
$X_{2}$ is $0\in\widetilde{K}(X_{2})$. As $X_{2}$ is a retract of $Y$ and
$X/A\cong Y/X_{2}$, the long exact sequence of $K$-theory for the pair
$\left(  Y,X_{2}\right)  $ yields the following short exact sequence:%
\[
0\rightarrow\widetilde{K}(X/A)\rightarrow\widetilde{K}(Y)\rightarrow
\widetilde{K}(X_{2})\rightarrow0\text{.}%
\]
Hence, $E_{1,2}-E_{2,2}$ is the image of a unique element in $\widetilde
{K}(X/A)$, which we name $[E_{1},E_{2};\sigma]$. It can be shown that this
element is invariant under variation of $\sigma$ within its own homotopy class
of bundle isomorphisms; see \cite{Atiyah K}.

\medskip

\textsc{Lemma 2.1. }\textit{In the above notation, we have the following
identities in }$\widetilde{K}(X/A)$\textit{:}

\begin{enumerate}
\item $[E_{1},E_{2};\sigma]+[E_{1}^{\prime},E_{2}^{\prime};\sigma^{\prime
}]=[E_{1}\oplus E_{1}^{\prime},E_{2}\oplus E_{2}^{\prime};\sigma\oplus
\sigma^{\prime}].$

\item $[E_{1},E_{2};\sigma]+[E_{2},E_{3};\rho]=[E_{1},E_{3};\rho\circ\sigma].$

\item $\left[  E_{2},E_{1};\sigma^{-1}\right]  =-\left[  E_{1},E_{2}%
;\sigma\right]  .$
\end{enumerate}

\medskip

\textsc{Proof: }Part 1 is a direct consequence of the definition.

For Part 2, note that $\sigma\oplus\rho:E_{1}|_{A}\oplus E_{2}|_{A}\rightarrow
E_{2}|_{A}\oplus E_{3}|_{A}$ is homotopic through bundle isomorphisms to%

\[
\varphi=%
\begin{pmatrix}
0 & -\operatorname{Id}\\
\rho\circ\sigma & 0
\end{pmatrix}
:E_{1}|_{A}\oplus E_{2}|_{A}\rightarrow E_{2}|_{A}\oplus E_{3}|_{A}%
\]
via
\[
t\mapsto%
\begin{pmatrix}
\operatorname{Id} & 0\\
t\rho & \operatorname{Id}%
\end{pmatrix}%
\begin{pmatrix}
\operatorname{Id} & -t\rho^{-1}\\
0 & \operatorname{Id}%
\end{pmatrix}%
\begin{pmatrix}
\operatorname{Id} & 0\\
t\rho & \operatorname{Id}%
\end{pmatrix}%
\begin{pmatrix}
\sigma & 0\\
0 & \rho
\end{pmatrix}
,\quad t\in\lbrack0,1].
\]
Thus%
\begin{align*}
\lbrack E_{1},E_{2};\sigma]+[E_{2},E_{3};\rho]  &  =[E_{1}\oplus E_{2}%
,E_{2}\oplus E_{3};\sigma\oplus\rho]\text{\quad(by Part 1)}\\
&  =[E_{1}\oplus E_{2},E_{2}\oplus E_{3};\varphi]\text{\quad(by }\varphi
\simeq(\sigma\oplus\rho)\text{)}\\
&  =[E_{1}\oplus E_{2},E_{3}\oplus E_{2};(\rho\circ\sigma)\oplus
(-\operatorname{Id})]\text{\quad(by the formula for }\varphi\text{)}\\
&  =[E_{1},E_{3};\rho\circ\sigma]+[E_{2},E_{2};(-\operatorname{Id}%
)]\text{\quad(by Part 1)}\\
&  =[E_{1},E_{3};\rho\circ\sigma]+[E_{2},E_{2};\operatorname{Id}]
\end{align*}
where the last equality follows from the fact that $\operatorname{Id}$ and
$-\operatorname{Id}$ are homotopic through isomorphisms of \textit{complex}
vector bundles. It is trivial that $[E_{2},E_{2};\operatorname{Id}]=0$ (which
can be checked in $\widetilde{K}(Y)$). Part 2 then follows.

Part 3 follows directly from Part 2.$\hfill\square$

\medskip

Given an isomorphism of complex vector bundles $\sigma:E_{1}\rightarrow E_{2}%
$, let $\sigma^{\ast}:E_{2}^{\ast}\rightarrow E_{1}^{\ast}$ denote the dual
isomorphism; if $E_{1}$ and $E_{2}$ are equipped with Hermitian metrics, let
$\widehat{\sigma}:E_{2}\rightarrow E_{1}$ denote the adjoint bundle map.
(Given Hermitian vector spaces $V$ and $W$ and a complex linear map
$g:V\rightarrow W$, the adjoint $\widehat{g}:W\rightarrow V$ is defined by
$\langle gv,w\rangle=\langle v,\widehat{g}w\rangle\ $for $v\in V$ and $w\in
W$. The adjoint of a bundle map can be fiberwise defined.) Since any complex
linear map $g:V\rightarrow W$ remains complex linear when viewed as a map
$g:\overline{V}\rightarrow\overline{W}$ where $\overline{V}$ and $\overline
{W}$ are the conjugates $V$ and $W$, we may view $\widehat{g}$ as $\widehat
{g}:\overline{W}\rightarrow\overline{V}$. In the same way, we may view
$\widehat{\sigma}$ as $\widehat{\sigma}:\overline{E_{2}}\rightarrow
\overline{E_{1}}$.

\medskip

\textsc{Lemma 2.2. }\textit{Suppose that }$E_{1}$\textit{\ and }$E_{2}%
$\textit{\ are complex vector bundles equipped with Hermitian metrics and that
}$\sigma:E_{1}|_{A}\rightarrow E_{2}|_{A}$\textit{\ is an isomorphism of
complex vector bundles. Then, in }$\widetilde{K}(X/A)$\textit{, we have:}

\begin{enumerate}
\item $[\overline{E_{2}},\overline{E_{1}};\widehat{\sigma}]=[E_{2}^{\ast
},E_{1}^{\ast};\sigma^{\ast}].$

\item $\left[  E_{1},E_{2};\sigma\right]  ^{\ast}=\left[  E_{1}^{\ast}%
,E_{2}^{\ast};\left(  \sigma^{\ast}\right)  ^{-1}\right]  .$

\item $\left[  E_{2},E_{1};\widehat{\sigma}\right]  =-\left[  E_{1}%
,E_{2};\sigma\right]  .$
\end{enumerate}

\medskip

\textsc{Proof: }For a complex linear map $g:V\rightarrow W$ between
finite-dimensional Hermitian vector spaces, we have the following commutative
diagram%
\[%
\begin{array}
[c]{ccc}%
\overline{W} & \overset{\widehat{g}}{\longrightarrow} & \overline{V}\\
\downarrow &  & \downarrow\\
W^{\ast} & \overset{g^{\ast}}{\longrightarrow} & V^{\ast}%
\end{array}
\]
where the vertical maps are the \textit{complex linear} isomorphisms%
\[
w\mapsto\left(  \langle-,w\rangle:W\rightarrow\mathbb{R}\right)  \in W^{\ast
}\text{ and }v\mapsto\left(  \langle-,v\rangle:V\rightarrow\mathbb{R}\right)
\in V^{\ast}.
\]
Part 1 then follows from the corresponding diagram of isomorphisms of bundles
over $A$.

For Part 2, recall that the bundle $E_{i,j}$ over $Y$ is given by gluing the
bundle $E_{i}$ over $X_{1}$ with $E_{j}$ over $X_{2}$ via a bundle isomorphism
$\varepsilon_{i,j}:E_{i}|_{A_{1}}\rightarrow E_{j}|_{A_{2}}$. It follows that
the bundle $E_{i,j}^{\ast}$ is formed by gluing $E_{i}^{\ast}$ with
$E_{j}^{\ast}$ via $\varepsilon_{i,j}^{\ast}:E_{j}^{\ast}|_{A_{2}}\rightarrow
E_{i}^{\ast}|_{A_{2}}$, or equivalently via $(\varepsilon_{i,j}^{\ast}%
)^{-1}:E_{i}^{\ast}|_{A_{1}}\rightarrow E_{j}^{\ast}|_{A_{2}}$. Part 2 readily follows.

Part 3 follows from a string of equalities:%
\begin{align*}
\left[  E_{2},E_{1};\widehat{\sigma}\right]   &  =[\left(  \overline{E_{2}%
}\right)  ^{\ast},\left(  \overline{E_{1}}\right)  ^{\ast};\sigma^{\ast}%
]\quad\text{(by Part 1)}\\
&  =\left[  \left(  \overline{E_{2}}\right)  ^{\ast},\left(  \overline{E_{1}%
}\right)  ^{\ast};\left(  \left(  \sigma^{-1}\right)  ^{\ast}\right)
^{-1}\right] \\
&  =\left[  \overline{E_{2}},\overline{E_{1}};\sigma^{-1}\right]  ^{\ast}%
\quad\text{(by Part 2)}\\
&  =-\left[  \overline{E_{1}},\overline{E_{2}};\sigma\right]  ^{\ast}%
\quad\text{(by Part 3 of Lemma 2.1)}\\
&  =-\left(  \overline{\left[  E_{1},E_{2};\sigma\right]  }\right)  ^{\ast}\\
&  =-\left[  E_{1},E_{2};\sigma\right]
\end{align*}
where the last equality follows from the bundle equivalence between the
conjugate of a complex vector bundle and its dual.$\hfill\square$

\bigskip

\begin{center}
\textsc{3. Proof of the Main Theorem}
\end{center}

\medskip

We now collect the preliminary results to prove our Main Theorem.

\medskip

\textsc{Proof of the Main Theorem: }Let $M$ be an even-dimensional closed
smooth manifold.

Because $\dim M$ is even, $\Lambda^{\text{even}}$ and $\Lambda^{\text{odd}}$
are both invariant under the involution $\tau$. Thus, $\tau$ decomposes
$\Lambda^{\text{even}}$ and $\Lambda^{\text{odd}}$ into ($\pm1$)-eigenbundles,
resulting in the following decomposition of $\Lambda^{\ast}$:%
\begin{align*}
\Lambda^{\ast}  &  =\Lambda^{\text{even}}\oplus\Lambda^{\text{odd}}\\
&  =(\Lambda^{\text{even},+}\oplus\Lambda^{\text{even},-})\oplus
(\Lambda^{\text{odd},-}\oplus\Lambda^{\text{odd},+})\text{.}%
\end{align*}
Likewise, $\Omega^{\text{even}}$ and $\Omega^{\text{odd}}$ can be decomposed
into ($\pm1$)-eigenspaces of $\tau$, resulting in the following decomposition
of $\Omega^{\ast}$:%
\begin{align}
\Omega^{\ast}  &  =\Omega^{\text{even}}\oplus\Omega^{\text{odd}}\nonumber\\
&  =(\Omega^{\text{even},+}\oplus\Omega^{\text{even},-})\oplus(\Omega
^{\text{odd},-}\oplus\Omega^{\text{odd},+})\text{.} \tag{1}%
\end{align}
Then, $D^{0}$ can be diagonalized as follows:%
\[
D^{0}=D^{0,+}\oplus D^{0,-}:\Omega^{\text{even},+}\oplus\Omega^{\text{even}%
,-}\rightarrow\Omega^{\text{odd},-}\oplus\Omega^{\text{odd},+}\text{.}%
\]

Regrouping the summands in (1), we have%
\begin{equation}
\Omega^{\ast}=(\Omega^{\text{even},+}\oplus\Omega^{\text{odd},+})\oplus
(\Omega^{\text{odd},-}\oplus\Omega^{\text{even},-})\text{.} \tag{2}%
\end{equation}
Clearly,%
\[
(\Omega^{\text{even},+}\oplus\Omega^{\text{odd},+})\subset\Omega^{+}%
\quad\text{and}\quad(\Omega^{\text{odd},-}\oplus\Omega^{\text{even},-}%
)\subset\Omega^{-}\text{ .}%
\]
Since $\Omega^{\ast}=\Omega^{+}\oplus\Omega^{-}$, the decomposition (2)
implies%
\[
\Omega^{+}=\Omega^{\text{even},+}\oplus\Omega^{\text{odd},+}\text{\quad
and\quad}\Omega^{-}=\Omega^{\text{odd},-}\oplus\Omega^{\text{even},-}\text{ .}%
\]
Then, $D^{+}$ can be diagonalized as follows:%
\[
D^{+}=D^{+,0}\oplus D^{+,1}:\Omega^{\text{even},+}\oplus\Omega^{\text{odd}%
,+}\rightarrow\Omega^{\text{odd},-}\oplus\Omega^{\text{even},-}\text{.}%
\]

Since $\Omega^{\text{even/odd},\pm}=\Gamma(\Lambda^{\text{even/odd},\pm})$, we
have%
\begin{align*}
\sigma(D^{0})  &  =\sigma(D^{0,+})\oplus\sigma(D^{0,-}):\pi^{\ast}%
\Lambda^{\text{even},+}\oplus\pi^{\ast}\Lambda^{\text{even},-}\rightarrow
\pi^{\ast}\Lambda^{\text{odd},-}\oplus\pi^{\ast}\Lambda^{\text{odd},+}\text{
,}\\
\text{and\quad}\sigma(D^{+})  &  =\sigma(D^{+,0})+\sigma(D^{+,1}):\pi^{\ast
}\Lambda^{\text{even},+}\oplus\pi^{\ast}\Lambda^{\text{odd},+}\rightarrow
\pi^{\ast}\Lambda^{\text{odd},-}\oplus\pi^{\ast}\Lambda^{\text{even},-}\text{
.}%
\end{align*}
By Part 1 of Lemma 2.1,%
\[
\left[  \sigma(D^{0})\right]  =\left[  \sigma(D^{0,+})\right]  +\left[
\sigma(D^{0,-})\right]  \quad\text{and}\quad\left[  \sigma(D^{+})\right]
=\left[  \sigma(D^{+,0})\right]  +\left[  \sigma(D^{+,1})\right]  \text{.}%
\]

It is a matter of definition that%
\[
D^{0,+}=D^{+,0}\text{\quad and\quad}D^{0,-}=\widehat{D^{+,1}}\text{.}%
\]
Thus,%
\[
\sigma(D^{0,+})=\sigma(D^{+,0})\text{\quad and\quad}\sigma(D^{0,-}%
)=\sigma(\widehat{D^{+,1}})=\widehat{\sigma(D^{+,1})}\text{.}%
\]
Therefore%
\[
\left[  \sigma(D^{0,+})\right]  =\left[  \sigma(D^{+,0})\right]
\]
and, by Part 3 of Lemma 2.2,%
\[
\left[  \sigma(D^{0,-})\right]  =\left[  \widehat{\sigma(D^{+,1})}\right]
=-\left[  \sigma(D^{+,1})\right]  \text{.}%
\]

In summary,%
\begin{align*}
\left[  \sigma(D^{0})\right]  +\left[  \sigma(D^{+})\right]   &  =\left[
\sigma(D^{0,+})\right]  +\left[  \sigma(D^{0,-})\right]  +\left[
\sigma(D^{+,0})\right]  +\left[  \sigma(D^{+,1})\right] \\
&  =2\left[  \sigma(D^{0,+})\right]
\end{align*}

$\hfill\square$

\bigskip

\begin{center}
\textsc{4. Discussion}

\medskip
\end{center}

By the $K$-theoretic Thom isomorphism theorem for complex vector bundles,
$K(T_{\mathbb{C}}^{\ast}M)$ is a free $K(M)$-module of rank 1 generated by the
$K$-theory Thom class. In this way, the Main Theorem admits interpretation in
$K(M)$. The $K$-theory Thom class of $T_{\mathbb{C}}^{\ast}M\rightarrow M$ is
simply the element%
\[
\lbrack\pi^{\ast}\Lambda^{\text{even}},\pi^{\ast}\Lambda^{\text{odd}};\mu]\in
K(T_{\mathbb{C}}^{\ast}M)
\]
with%
\[
\mu(v,\omega)=v\wedge\omega-\iota_{v}\omega\quad\text{for }v\in(T_{\mathbb{C}%
}^{\ast}M)_{x}\text{ and }\omega\in(\Lambda^{\text{even}})_{x}%
\]
where $\iota_{v}:\Lambda^{k}\rightarrow\Lambda^{k-1}$ is the interior
multiplication by $v$. This is precisely the $K$-theoretic Euler class, i.e.,
the symbol class of the de Rham operator $\left[  \sigma(D^{0})\right]  \in
K(T_{\mathbb{C}}^{\ast}M)$, which corresponds (via Thom isomorphism) to $1\in
K(M)$. Thus, letting $S\in K(M)$ denote the element corresponding (via the
Thom isomorphism) to $\left[  \sigma(D^{+})\right]  \in K(T_{\mathbb{C}}%
^{\ast}M)$ and translating the Main Theorem into a statement in $K(M)$, we
have that $S=1+2x$ for some $x\in K(M)$.

\smallskip

When $M$ admits an orientation-preserving smooth action by a compact Lie group
$G$, the Main Theorem admits an equivariant generalization. (This is due to a
suggestion by Sylvain Cappell.) Equip $M$ with a $G$-invariant metric, i.e., a
metric with respect to which $G$ acts isometrically. Then, $\left[
\sigma(D^{0})\right]  $ and $\left[  \sigma(D^{+})\right]  $ can both be
interpreted as $G$-equivariant $K$-theory classes, i.e., elements of
$K_{G}(T_{\mathbb{C}}^{\ast}M)$, and the proof for the Main Theorem can be
adapted to show the following.

\medskip

\textsc{Theorem. }\textit{Suppose that }$\dim M$\textit{ is even and that a
compact Lie group }$G$\textit{ acts smoothly on }$M$\textit{ preserving
orientation. Then }$\left[  \sigma(D^{0})\right]  \equiv\left[  \sigma
(D^{+})\right]  \operatorname{mod}2K_{G}(T_{\mathbb{C}}^{\ast}M)$.

\medskip

When $G$ is finite, we may apply $G$-$\operatorname*{Index}:K_{G}%
(T_{\mathbb{C}}^{\ast}M)\rightarrow R(G)$ to deduce from this result the mod-2
equivalence of the equivariant Euler characteristic and multi-signature (as
virtual complex representations of $G$).

\smallskip

In closing, we mention a few related works. One is \cite{R}, which shows that
the $K$-homology class of the de Rham operator is trivial; another is
\cite{LR}, which discusses the equivariant $KO$-theoretic Euler
characteristic. A more recent work, \cite{RW}, shows among other results that
the mod-8 reduction of the $K$-homology class of the signature operator is an
oriented homotopy invariant.

\bigskip

\noindent%
\begin{tabular}
[c]{lll}%
Department of Mathematics & \qquad\qquad\qquad & (In Transition)\\
Indiana University &  & \\
Bloomington, IN 47405 &  & \\
United States of America &  & United States of America\\
Email:\quad jfdavis@indiana.edu &  & Email:\quad pd260@nyu.edu
\end{tabular}

\end{document}